\documentclass[10pt,reqno]{amsart} 

\usepackage{amssymb,latexsym}
\usepackage{cite} 

\usepackage[height=190mm,width=130mm]{geometry} 

\theoremstyle{plain}
\newtheorem{theorem}{Theorem}
\newtheorem{lemma}{Lemma}

\theoremstyle{definition}

\theoremstyle{remark}

\newtheorem{example}{\bf Example}

\newcommand{\Z}{\mathbb{Z}}
\newcommand{\C}{\mathbb{C}}

\numberwithin{equation}{section} 

\newcommand{\rmap}{\longrightarrow}

\newcommand{\g}{\ensuremath{\Gamma}}

\newcommand{\ps}{{\raise 1pt\hbox{\tiny (}}}

\newcommand{\pss}{{\raise 1pt\hbox{\tiny [}}}
\newcommand{\pdd}{{\raise 1pt\hbox{\tiny ]}}}
\newcommand{\pd}{{\raise 1pt\hbox{\tiny )}}}

\newcommand{\bs}{{\raise 1pt\hbox{\tiny [}}}
\newcommand{\bd}{{\raise 1pt\hbox{\tiny ]}}}

\def\cross{\mathinner{\mathrel{\raise0.8pt\hbox{$\scriptstyle>$}}
                 \joinrel\mathrel\triangleleft}}

\usepackage{citehack}
\usepackage{amsmath,amsthm}
\usepackage{amsfonts}
\usepackage{amssymb}
\usepackage{longtable}
\usepackage[matrix,arrow,curve]{xy}

\usepackage{lipsum}

\usepackage{amsfonts}
\usepackage{amsmath}
\usepackage{amscd}
\usepackage{amssymb}
\usepackage{amsthm}
\usepackage{bbm}
\usepackage{bm}
\usepackage{enumerate}
\usepackage{enumitem}
\usepackage{hyperref}
\usepackage{latexsym}
\usepackage{mathabx}
\usepackage{mathdots}
\usepackage{mathtools}

\usepackage{amsmath}
\usepackage{amssymb}
\usepackage{amsthm}
\usepackage{bbm}
\usepackage{enumerate}
\usepackage{latexsym}
\usepackage{mathdots}
\usepackage{geometry}

\usepackage{hyperref}  
 
\usepackage{amsmath}
\usepackage{amssymb}
\usepackage{amsthm}
\usepackage{bbm}
\usepackage{enumerate}
\usepackage{latexsym}
\usepackage{mathdots}
\usepackage{geometry}

\def\K{\mathcal{K}}

\usepackage{stackrel}

\newcommand{\be}{\begin{equation}}
\newcommand{\ee}{\end{equation}}

\newcommand{\nn}{\nonumber \\}

\newcommand{\wt}{\mbox{\rm wt}\ }

\newcommand{\nc}{\newcommand}
\nc{\cali}{\mathcal}
\nc{\on}{\operatorname}
\nc{\Wick}{{\mb :}}

\nc{\ddz}{\frac{\partial}{\partial z}}
\nc{\ch}{\mbox{ch}}
\nc{\Oo}{{\cali O}}
\nc{\cond}{|\,}
\nc{\bib}{\bibitem}
\nc{\pone}{\Pro^1}
\nc{\pa}{\partial}
\nc{\arr}{\rightarrow}
\nc{\larr}{\longrightarrow}
\nc{\ket}{\rangle}
\nc{\bra}{\langle}
\nc{\gam}{\bar{\gamma}}
\nc{\q}{\widetilde{Q}}
\nc{\ep}{\epsilon}
\nc{\su}{\widehat{{\mf s}{\mf l}}_2}
\nc{\sw}{{\mf s}{\mf l}}
\nc{\h}{{\mf h}}
\nc{\n}{{\mf n}}
\nc{\ab}{\mf{a}}
\nc{\is}{{\mb i}}
\nc{\js}{{\mb j}}
\nc{\bi}{\bibitem}
\nc{\He}{{\cali H}}
\nc{\inv}{^{-1}}
\nc{\ol}{\overline}
\nc{\wh}{\widehat}
\nc{\dst}{\displaystyle}

\nc{\delt}{\partial_t}
\nc{\ddt}{\frac{\partial}{\partial t}}
\nc{\delx}{\partial_x}
\nc{\mb}{\mathbf}
\nc{\mf}{\mathfrak}

\nc{\mbb}{\mathbb}
\nc{\Ctt}{\C((t))}
\nc{\Ct}{\C[t,t\inv]}

\nc{\ghat}{\wh{\g}}

\nc{\un}{\underline}
\nc{\mc}{\mathcal}
\nc{\BB}{{\mc B}}
\nc{\bb}{{\mf b}}
\nc{\kk}{{\mf k}}
\nc{\frob}{\times}
\nc{\sm}{\setminus}
\nc{\Pp}{{\mathbb P}^1}
\nc{\Aa}{{\mc A}}

\nc{\AutO}{\on{Aut}\Oo}
\nc{\AUTO}{\un{\on{Aut}}\Oo}
\nc{\AUTK}{\un{\on{Aut}}\K}
\nc{\Heout}{\He_{\out}}
\nc{\Hetil}{{\widetilde\He}}
\nc{\wb}{\overline}

\nc{\Res}{\on{Res}}
\nc{\pitil}{\Pi}
\nc{\Ctil}{\wt{C}}
\nc{\auto}{\on{Aut} \Oo}
\nc{\phitil}{\wt{\phi}}
\nc{\gz}{\g_{\vec z}}
\nc{\tensorM}{\bigotimes_{i=1}^N{\mathbb M}_i}
\nc{\tensorW}{\bigotimes_{i=1}^N W_{\nu_i,k}}
\nc{\out}{\on{out}}

\nc{\m}{{\mathfrak m}}

\nc{\gx}{\g^0_{\vec x}}

\nc{\hx}{\He^0_{\vec x}}
\nc{\tensorpi}{\pi_{\nu_1,\ldots,\nu_N}^\kappa}
\nc{\Phizw}{\Phi_{\vec w}({\vec z})}
\nc{\Pro}{{\mathbb P}}

\nc{\De}{\Delta}

\nc{\us}{\underset}

\nc{\Ll}{\mc L}
\nc{\dR}{\on{dR}}

\nc{\T}{{\mc T}}

\nc{\Xn}{\overset{\circ}X{}^n} \nc{\Dn}{\overset{\circ}D{}^n}
\nc{\Dxn}{\overset{\circ}D{}^n_x} \nc{\varphitil}{\wt{\varphi}}

\nc{\lf}{{\mf l}}
\nc{\GL}{{}^L G}
\nc{\Vir}{\on{Vir}}

\begin{document}
\title[Invariant classes for families of complexes]  
{Invariant classes for families of complexes}  

\author{Daniel Levin${}^\&$}
\address{ ${}^\&$ Mathematical Institute \\ University of Oxford \\ Andrew Wiles Building \\
 Radcliffe Observatory Quarter (550) \\ Woodstock Road \\ Oxford \\ OX2 6GG \\ United Kingdom
}

\author{Alexander Zuevsky${}^*$}
\address{ ${}^*$ Institute of Mathematics \\ Czech Academy of Sciences\\ Zitna 25, 11567 \\ Prague\\ Czech Republic}

\email{levindanie@gmail.com}
\email{zuevsky@yahoo.com}

\begin{abstract}
We consider families of chain-cochain infinite complexes 
$\mathcal C$  
of spaces with elements depending on a number of parameters, and 
 endowed with 
a converging associative multiple product. 
The existence of left/right local/non-local square-vanishing ideals 
is assumed for subspaces of $\mathcal C$-spaces. 
We show that a set of differential and orthogonality relations 
together with coherence conditions
 on indices of a chain-cochain complex $\mathcal C$ elements   
generates families of graded differential algebras. 
With the appropriate orthogonality conditions on completions 
of $\mathcal C$ elements in the multiple product, 
 we define the equivalence classes of cohomology invariants.  
AMS Classification: 53C12, 57R20, 17B69 
\end{abstract}

\keywords{Cohomology invariants; 
 orthogonality condition for chain-cochain complexes; 
 graded differential algebras} 
\vskip12pt  

\maketitle

\begin{center}
{Conflict of Interest and Data availability Statements:}
\end{center}

The authors state that: 

1.) The paper does not contain any potential conflicts of interests. 

2.) The paper does not use any datasets. No dataset were generated during and/or analysed 
during the current study. 

3.) The paper includes all data generated or analysed during this study. 

4.) Data sharing is not applicable to this article as no datasets were generated or analysed during the current study.

5.) The data of the paper can be shared openly.  

6.) No AI is used to write this paper. 
\section{Introduction} 
\label{valued}
The cohomology of complexes endowed with associative products 
is an important way to study various algebraic and geometric structures. 
Computation of cohomology classes of invariants associated to 
such complexes is difficult problem. 
For differential forms considered on smooth manifolds,  
the Frobenius theorem for a distribution leads to the orthogonality condition.   
Motivated by the notion of integrability for differential forms on foliated manifolds 
\cite{Ghys, Ko},    
 we introduce in this paper the notion of orthogonality 
with respect to a multiple product of elements of a chain-cochain complex. 
 Differential and orthogonality conditions then 
endow a complex with the structure of 
a graded differential algebra. 
Assuming existence of square-vanishing ideals in subspaces of families of complexes,  
and requiring natural orthogonality conditions for completions of elements, 
cohomology classes
are derived and their independence is proven. 
The classes obtained constitute generalization of Godbillon-Vey classes 
for foliations \cite{Ghys, Ko}. 

As for possible applications of the material presented in this paper,
 we would like to mention 
computations of higher cohomology 
of grading-restricted vertex algebras \cite{Huang} 
endowed with a sewing procedure \cite{G, Yamada},    
constructions of more complicated 
cohomology invariants, and applications in differential geometry and algebraic topology.  
In particular, it will be possible to compute cohomology invariants of 
 vertex operator algebra bundles \cite{BZF}.  
It would be interesting to study possible applications of invariants we constructed
 to cohomology of manifolds, continual Lie algebras \cite{sv2},  
conformal field theory \cite{FMS, TUY}.  
In differential geometry there exist various
 approaches to the construction of cohomology classes 
 (cf., in particular, \cite{Losik}). 
The results proven in this paper 
 are useful in computations of foliations cohomology \cite{BG, BGG}. 
\section{Setup and the main results}
\label{setup}
Introduce a system of families of horizontal chain-cochain complexes   
$\mathcal C$ $=$ 
 $\left(C^{n, \kappa}_m \right.$, $\delta^{n, \kappa}_m$,
 $\left. \Delta^{n, \kappa}_m \right)$
 with $n$, $m$, $\kappa \in \Z$.   
The horizontal complexes with differentials
 $\delta^{n, \kappa}_m: C^{n, \kappa}_m \to C^{n+1, \kappa}_{m-1}$  
are supposed to be chain-cochain, 
i.e.,  
$\delta^{n+1, \kappa}_{m-1} \circ \delta^{n, \kappa}_m=0$, for all $\kappa \in \Z$.  
In certain cases we take the vertical complexes with differentials 
$\Delta^{n, \kappa}_m: C^{n, \kappa}_m \to C^{n, \kappa+1}_m$,  
 to be chain complexes with 
$\Delta^{n, \kappa}_m$ commuting with $\delta^{n', \kappa'}_{m'}$. 
When we write $c\in \mathcal C$ that means that an element $c$ belongs to 
a space of $\mathcal C$.  
Instead of families of $C^{n, \kappa}_m$-complexes for various $\kappa$, 
one can consider families $C^{\bm n}_{\bm m}$ with multiple indices 
  ${\bm n}=(n_1, \ldots, n_{M_1})$, ${\bm m}=(n_1, \ldots, n_{M_2})$,  
for some $M_1$ and $M_2 \ge 1$.  
Each space of the multiple complex has two types of parameter 
with respect to the action of differentials: 
one set of growing parameters $\bm n$, and   
one set of lowering parameters $\bm m$. 
 Such a construction will be considered in a separate paper. 

Let us give an example of a $\mathcal C$-complex with three parameters 
$n$, $m$, $\kappa \in \Z$ which corresponds  
to the diagram 
\begin{eqnarray}
\label{buzovaish}
&&  \qquad \qquad \quad   \vdots \qquad \qquad \qquad \qquad  \vdots
\nn
&& \qquad \qquad  \quad \downarrow  \Delta^{n, \kappa-1}_m
\qquad \; \qquad  \downarrow  \Delta^{n, \kappa-1}_{m-1}
\nn
&& \ldots \stackrel{\delta^{n-1, \kappa}_{m+1}} {\longrightarrow} \quad C^{n, \kappa}_m  
  \qquad \stackrel{\delta^{n, \kappa}_m} {\longrightarrow} \qquad  C^{n+1, \kappa}_{m-1}    
  \stackrel{\delta^{n+1, \kappa}_{m-1}} {\rmap} 
 \cdots  
\nn
&& \qquad \qquad \quad  \downarrow  \Delta^{n, \kappa}_m 
\qquad \quad   \qquad  \downarrow  \Delta^{n+1, \kappa}_{m-1}
\nn
&& \ldots \stackrel{\delta^{n-1, \kappa+1}_{m+1}} {\longrightarrow} \quad C^{n,\kappa+1}_m 
   \stackrel{\delta^{n, \kappa+1}_m} {\longrightarrow}  \quad C^{n+1,\kappa+1}_{m-1}    
  \stackrel{\delta^{n+1, \kappa+1}_{m-1}} {\rmap} 
 \cdots  
\nn
&& \qquad \qquad \quad    \downarrow  \Delta^{n, \kappa+1}_m
\qquad  \quad \quad  \downarrow  \Delta^{n+1, \kappa+1}_{m-1} 
\nn
&&  \qquad \qquad \quad   \vdots \qquad \qquad \qquad  \qquad  \vdots 
\end{eqnarray} 
Since the differentials $\delta^{n, \kappa}_m$ 
and $\Delta^{n, \kappa}_m$ are supposed to act on elements of 
appropriate spaces $C^{n, \kappa}_m$ with the same set of indices    
we will often omit indices and write $\delta$ and $\Delta$. 

For all combinations of $l$ elements of various spaces of $\mathcal C$ 
we assume existence of a formal multiple  
 associative product 
\begin{equation}
\label{tugo}
\cdot_l=\cdot_{i=1}^l: \times_{i=1}^l C^{n_i, \kappa_i}_{m_i} 
\to C^{n',\kappa'}_{m'},   
\end{equation}  
for $n'=\sum_{i=1}^l n_i$, $\kappa' =\sum_{i=1}^l \kappa_i$, 
$m'=\sum_{i=1}^l m_i$.  
It is assumed that a result of the $\cdot_l$-product
which is an element of a $\mathcal C$-space, 
is well-defined (e.g., converging) for all values of 
parameters defining the product.  
The associativity of a $\cdot_l$-product means  
$(a\cdot_l \ldots \cdot_{l'} b) \cdot_{l''} \ldots \cdot_{l'''}c
=a \cdot_l \ldots \cdot_{l'} (b \cdot_{l''} \ldots \cdot_{l'''} c)$. 
 for elements $a$, $b$, $c \in \mathcal C$.  

In certain cases, a product \eqref{tugo} may lead to coincidence of 
parameters of multiplied elements of a complex \eqref{buzovaish} spaces. 
Resulting product may not allow such a coincidence resulting in 
possible overcounting of the number of parameters  
(e.g., for a product of vertex operator algebra complexes \cite{Huang}). 
In order to avoid such a possibility, 
we take into account one coinciding elements/parameters only in the resulting 
product.  
For an element $\phi_i \in C^{n_i, \kappa_i}_{m_i}$ 
in a particular $\cdot_l$-product  
let $r_i$ and $t_i$ be the numbers of common parameters with other elements 
in the product
corresponding to upper and lower indices for $\phi_i$. 
In that case, the conditions for indices for the resulting element 
of a $\cdot_l$-product are 
$n_i = \sum_{j=1}^l n_j-r_j$, and $m_i = \sum_{j=1}^l n_j-t_j$. 

Note that $l$ can be infinite if we assure that the result of the product 
$\times_{i \ge 1}$ is converging. 
We call the product formal since elements of a complex $\mathcal C$
spaces can be formal geometric objects. 
Then the result of a $\cdot_l$-product is 
 a superposition 
of geometrical objects, e.g., Riemann surfaces.  
In that case, convergence of the product means that 
corresponding superposition leads to a well-defined formal object. 
Note that we do not assume commutation relations for elements of 
a $\mathcal C$-spaces to be known. 
All constructions of this paper are independent of actual commutation 
formulas for elements inside a $\cdot_l$-product. 
 We assume that   
in the spaces of a complex $\mathcal C$ 
there exist a net of subspaces constituting 
  (left, right) square-vanishing 
 ideals $\mathcal I(2) \subset \mathcal C$. 
I.e., the if for elements $a'$, $a \in \mathcal I(2) \subset \mathcal C$, 
satisfying the condition $\cdot_l\left(a', b, a\right)=0$, 
then there exists $c$, $c'$, $c''$, $d$, $d' \in \mathcal C$,   
$b=\cdot_l(a', c, d)$ or $b=\cdot_l(d', c', a)$, or 
$b=\cdot_l(a', c'', a)$, 
for $b \in \mathcal I_L$ or $b \in \mathcal I_R$ or $b \in \mathcal I_{L, R}$ 
correspondingly. 
We define (left, right) local $\mathcal I_l(2)$ and non-local 
$\mathcal I_{nl}(2)$ ideals of order $2$ as  
the set of $\theta$, $\theta' \in \mathcal C$ such that 
$\cdot_l(\ldots, \theta, \theta', \ldots )=0$ or
$\cdot_l(\ldots, \theta, \ldots, \theta', \ldots )=0$ correspondingly. 
The cases of higher order ideals $\mathcal I(k)$, $k >2$
 will be considered in another paper.  

Introduce the following notations.
Since differentials and elements of $\mathcal C$-complex 
are always assumed coherent, let us denote  
$\delta =\delta^{i,k}_m$,  $\Delta=\Delta^{i,k}_m$. 
  Let $d=d^{i,k}_m$ be a choice of $\delta$ or $\Delta$,  
and $\overline{d}$ denotes the alternative choice.  
In addition to all that above it is set that a product \eqref{tugo} 
is coherent with respect to all vertical and horizontal differentials 
of the diagram \eqref{buzovaish}. 
That means that for every $\kappa=\sum_{i=1}^l \kappa_i$,
 $n=\sum_{i=1}^l n_i$,  
 and $m=\sum_{i=1}^l m_i$, 
and corresponding $\phi \in C^{n, \kappa}_m$
 given by a $\cdot_l$-product, 
 the action of corresponding differential $d=d^{n, \kappa}_m$  
has Leibniz formula
\begin{equation}
\label{leibniz}
d^{n, \kappa}_m \phi 
= d^{n, \kappa}_m \cdot_l 
\left(\phi_1, \ldots, \phi_i, \ldots, \phi_l \right)=  
\sum\limits_{i=1}^l \cdot_l 
\left(\ldots, d^{n_i, \kappa_i}_{m_i}   
\phi_i, \ldots \right),   
\end{equation}
 with respect to the elements 
$d^{n_i, \kappa_i}_{m_i} \phi_i$,   $\phi_i \in C^{n_i, \kappa_i}_{m_i}$.  

For an element $\gamma \in \mathcal C$ and a set     
$J=(j_1, \ldots j_l)$ of choices $j_i$ of $(0, 1)$,  
let us require that for $l$ chain-cochain spaces  
$C^{n_i, \kappa_i}_{m_i}$ of the complex \eqref{buzovaish},   
 there exist subspaces 
$C^{'n_i, \kappa_i}_{m_i} \subset C^{n_i, \kappa_i}_{m_i}$    
 such that for all $\phi_i \in C^{'n_i, \kappa_i}_{m_i}$,   
\begin{equation}
\label{purgo}
 \delta \gamma
=\cdot_l \left(\delta^{j_1} \phi_1, \ldots, \delta^{j_i} \phi_i, \ldots, 
\delta^{j_l} \phi_l, \right).
\end{equation} 
 We call \eqref{purgo} a differential condition. 
Symmetrizing with respect to all possible choices of
 $J_k$, $k \ge 1$, with a fixed $l$,    
we obtain a set of differential conditions 
\begin{equation}
\label{mordo}
\left\{  Symm_{J_{l_k}} \left\{ \delta \gamma_k=
 \; \cdot_{l_k} \left(\delta^{j_{1, k}} \phi_{1, k}, \ldots, \delta^{j_{i, k}} 
\phi_{i, k}, \ldots, \delta^{j_{l_k, k}} \phi_{l_k, k}\right) \right\} \right\}, 
\end{equation}
for $\gamma_k \in \mathcal C$, $l_k \ge 1$, $k \ge 1$.  
With all $\delta \gamma_k=0$,    
 we call \eqref{mordo} the set of orthogonality conditions. 
For each differential condition above of the set \eqref{mordo} 
  we have the coherence condition for 
corresponding upper and lower indices of $\mathcal C$-spaces, i.e., 
$n_{\delta \gamma_k}=\sum_{i=1}^{l_k} n_{\delta^{j_i, k}\phi_{i, k}}
- r_{\delta^{j_i, k}\phi_{i, k}}$,  
$m_{\delta \gamma_k}=\sum_{i=1}^{l_k} m_{\delta^{j_i, k}\phi_{i, k}}
- t_{\delta^{j_i, k}\phi_{i, k}}$, 
and 
$\kappa_{\delta \gamma_k}=\sum_{i=1}^{l_k} \kappa_{\delta^{j_i, k}\phi_{i, k}}$. 
We will skip these coherence relations 
for all differential conditions below. 

The first result of this paper is the following 
Lemma \ref{kuzya} which will be illustrated in Section \ref{graded}. 
\begin{lemma}
\label{kuzya} 
Actions of the differentials $d$, $\overline{d}$ to 
 differential and orthogonality conditions  
with respect to a $\cdot_l$-multiplication \eqref{tugo}, $l \ge 0$, 
 endow families of 
 chain-cochain complexes  
$\mathcal C$ 
with the structure of a multiple graded 
 differential algebra.  
\end{lemma} 
For a set of $p$ elements $\phi_i$, $1\le i \le p$,  
in a $\cdot_l$-product 
\begin{equation}
\label{mama}
\cdot_l(\Phi_1, \phi_1, \Phi_2, \ldots, 
\Phi_i, \phi_i, \Phi_{i+1}, \ldots, \Phi_p, \phi_p, \Phi_{p+1}),  
\end{equation}
we call the union of sets of $\mathcal C$-elements
 $\left\{ \Phi_j \right\}$, $1 \le j\le p+1$, 
the completion of a set of elements $\phi_i$ in a $\cdot_l$-product. 
Due to associativity of the $\cdot_l$-product
 we can consider each set $\Phi_j$ as one element.
The completion $\left\{ \Phi_j \right\}$, $1 \le j\le p+1$, 
in \eqref{mama} for a set of elements $\phi_i$, $1 \le i \le p$,  
with respect to the differential $d$  
 is called closed if 
\begin{equation}
\label{papa}
\sum\limits_{j=1}^{p+1} d_j.\cdot_l(\Phi_1, \phi_1, \Phi_2, \ldots, 
\Phi_i, \phi_i, \Phi_{i+1}, \ldots, \Phi_p, \phi_p, \Phi_{p+1})=0,  
\end{equation}
where $d_j$ acts only on $\Phi_j$-elements. 
Denote by $\xi_{\sigma(j)}$, $1\le j \le 3$
the permutations with repetitions of the set of two $\mathcal C$-element $(\phi, \psi)$.  
Let us define an auxiliary operator 
$\epsilon: \mathcal C\times \mathcal C \to \mathcal C$, 
such that $\epsilon(.,.)=(.,.)$, when $d$ or $\overline{d}$ is chain-cochain
and $\epsilon(.,.)={\rm Id}_{\mathcal C}$ otherwise.  
Consider the case when $\overline{d}d \phi\ne 0$, for all $\phi \in \mathcal C$.  
We then formulate then the second result of this paper: 
\begin{theorem} 
\label{marcela} 
For arbitrary elements $\phi$, $\psi \in \mathcal C$, 
 such that $d\phi$, $d \psi \in \mathcal I_{nl}(2)$,  
any set of closed completions $\Phi_j \in \mathcal C$, $1 \le j \le 4$, 
with respect to the choice $d$ of differentials $\delta$ or $\Delta$  
for the permutations $\xi_{\sigma(1\le j_1, j_2, j_3 \le 3)}$  
with repetitions of the set $(\phi, \psi)$, i.e., 
\begin{equation}
\label{tuba}
  \left\{ { Symm \atop {\small 1 \le j_1, j_2, j_3 \le 3}} \left\{
0=\sum\limits_{j=1}^4 d_j.\cdot_l \left( \epsilon\left(\Phi_1, 
 \overline{d_i}   
\xi_{\sigma(j_1)}\right),  
 \Phi_2, d_i \xi_{\sigma(j_2)},   
\Phi_3, \xi_{\sigma(j_3)},  
\Phi_4\right) \right\} \right\},  
\end{equation}
the cohomology classes are given by  
\begin{equation}
\label{kordo}
 Symm_J \;\left\{ \;  \left[ \cdot_l 
\left(\epsilon\left(\Phi_1, \overline{d} d \phi \right),     
  \Phi_2, d \phi, 
\Phi_3,  \phi, \Phi_4  
 \right) 
\right]\; \right\},    
\end{equation}
 independent on the choice of 
 $\phi \in \mathcal C$. 
\end{theorem}
The proof of Theorem \ref{marcela} is given in Section \ref{invariants}. 
The result of Theorem  \ref{marcela} generalizes results of \cite{Ghys}
for Godbillon-Vey invariants for codimension one foliations, 
and results of \cite{Ko}
for differential forms defined on foliations. 
The set of conditions \eqref{tuba} generalizes 
corresponding orthogonality conditions assumed in \cite{Ghys, Ko}. 
Due to the associativity property of the $\cdot_l$-product mention above, 
we will skip the notation $\cdot_l$ and denote all products for various $l$ 
as $(\ldots, \ldots, \ldots)$. 
Since all actions of differentials of are supposed to be coherent 
with the indexing of $\mathcal C$-spaces, we will skip also the upper and lower 
indices from $\delta^{n, \kappa}_m$. 
Let us underline that we work with a complex spaces as 
grading subspaces of an associative algebra. 
No commutation relations for element of various subspaces are assumed known. 
In particular, the multiple product $\cdot_l$ can be introduces for 
geometric objects of general kind. 
Thus, we keep the constructions in this paper independent of commutation relations. 
\section{Graded differential algebras associated to families of chain-cochain 
complexes} 
\label{graded} 
The way of manipulation with sequences 
of differential and orthogonality conditions 
 follows from 
the chain-cochain and Leibniz \eqref{leibniz} conditions for differentials, as well as 
the assumption that some of $\mathcal C$ belong to the ideals 
$\mathcal I(2)$ or $\mathcal I_{nl}(2)$.  
Now we provide a specific example illustrating 
  Lemma \ref{kuzya}. 
\begin{example} 
Due to associativity,  
we consider each differential in a $\cdot_l$ product    
 for $1 \le i \le l$ separately. 
A particular element $\delta^{j_i} \phi_i$ may be situated 
at an arbitrary position in the $\cdot_l$-product in 
a differential or orthogonality condition defined by 
 a set $J=(j_1, \ldots, j_l)$, where $j_i=(0, 1)$ is a choice of $0$ of $1$. 
In this example we consider splittings of the product into two and three parts. 
In general, starting from a particular orthogonality or differential condition, 
we act consequently by $d$ and $\overline{d}$ differentials.  
For elements belonging 
to the left of right ideals $\mathcal I_l(2)$ 
or $\mathcal I_{nl}(2)$ of order two,  
the squares of any pair of their elements vanish. 
 Using the $\mathcal I(2)$ ideal vanishing properties
applied to orthogonality conditions, 
we express particular elements $\delta^{j_i} \phi_i$ 
in terms of other $\mathcal C$ elements. 
Finally, by taking into account corresponding coherence conditions,  
we arrive to the full structure of differential relations. 
A sequence of relations does not stop 
as long as coherence conditions on indices are fulfilled, 
or until the sequence gives identical zero. 

Let us reproduce the general structure of relations 
following from particular orthogonality conditions.  
It is quite clear that differential conditions 
always follow from some orthogonality conditions. 
In this example we consider only two first types of cases where 
differential conditions with respect to one differential $\delta$ 
 follow
 from its applications 
to the orthogonality conditions 
of the form $(\Phi, \phi)$ and $(\Phi', \phi, \Phi)$ with 
one picked element $\phi \in \mathcal C$ in a $\cdot_l$-product. 
The consideration below extends to the case of higher number 
of picked elements in a product. 
The case of arbitrary action of both differentials $d$ and $\overline{d}$
on arbitrary number of elements in a multiple product will be considered 
in a separate paper. 

Let us start with the form $(00)$ and the differential $\delta$. 
Similar computations can be done for the differential $\Delta$ 
and for combinations of these two differentials. 
The differential part: 
\begin{eqnarray*}
&& (00): \; 0=(\Phi', \Phi); \; 
  0=(\delta\Phi', \Phi) + (\Phi', \delta\Phi); \; 
 0=(\delta\Phi', \delta \Phi)
 + (\delta\Phi', \delta\Phi); \; 
 0 \;  \mbox{or} \; (II). 
\end{eqnarray*}
The square vanishing-part
\begin{eqnarray*}
&& (00): \; 0=(\Phi', \Phi); \; \Phi'=(\beta', \Phi); 
\\
&& 
1) \; \delta\Phi'=  (\delta\beta', \Phi) 
+ (\beta', \delta \Phi); \; 0 \;  \mbox{or} \; (II); \;   
2) \; \delta\Phi = (\Phi', \beta); \; 
(II). 
\end{eqnarray*}
The case $(I0)$, differential part: $0=(\delta \phi, \Phi); \; (II)$.    
Non-differential part: 
\begin{eqnarray*}
1)\; \delta \phi =(\beta', \Phi); \; 
 0=(\delta \beta', \Phi) + (\beta', \delta \Phi); \;  (II); \; 
2) \; \Phi=(\delta \phi, \beta); \; 
 \delta \Phi=(\delta \phi, \delta \beta); \; 0. 
\end{eqnarray*}
Case $(0I)$: differential part: 
$0=(\Phi', \delta \phi); \; (II)$.  
Non-differential part: 
\begin{eqnarray*}
 1)\; \Phi'=(\beta', \delta \phi); \; \delta \Phi'=(\delta \beta', \delta \phi); \; 0; \;  
 2)\; \delta \phi=(\Phi', \beta); \;  
0=(\delta \Phi', \beta)+ (\Phi', \delta \beta), \; (II). 
\end{eqnarray*}
Let us now consider the case $(II) \to (II)$. 
The differential part is trivial. 
The square vanishing-part is: 
\begin{eqnarray*}
&& (II): \; 0=(\delta \Phi', \delta \Phi); \; 
  1)\; \; \delta \Phi'=(\beta', \delta \Phi); \; 
  (II): \; 0=(\delta \beta', \delta \Phi); \; 
 2)\; \delta \Phi= (\delta \Phi', \beta); \; (II). 
\end{eqnarray*}
Let us continue with the differentiation of the form $(000)$: 
\begin{eqnarray*}
&& (000): \; 0=(\Phi', \phi, \Phi), \; 
 0= (\delta\Phi', \phi, \Phi) + (\Phi', \delta \phi,  \Phi) + (\Phi', \phi, \delta \Phi),  
\\
&&
0= (\delta\Phi', \delta\phi, \Phi) + (\delta\Phi', \phi, \delta \Phi) 
 + (\delta \Phi', \delta \phi, \Phi) +  (\Phi', \delta \phi, \delta \Phi)
+ (\delta \Phi', \phi, \delta \Phi) + (\delta \Phi', \delta \phi, \delta \Phi),  
\\
&&
(III): \; 0= (\delta\Phi', \delta\phi, \delta \Phi); \; 
1) \; \delta\Phi'= (\beta', \delta\phi); \; (II): \;  0= (\delta \beta', \delta\phi);   
\\
&&
2a)\; \delta\phi = (\delta\Phi', \gamma); \; (II): 0=(\delta\Phi', \delta \gamma); \; 
2b) \; \delta\phi = (\gamma', \delta \Phi); \; (II): \; 0= (\delta \gamma', \delta \Phi), 
\\
&& 3)\;  \delta \Phi = (\beta,  \delta\phi ); \; (II): \; 0 = (\delta \beta,  \delta\phi ). 
\end{eqnarray*}
Next we continue with the square vanishing-part of elements of $(000)$
and differentiations:  
\begin{eqnarray*}
&& (000): \; 0=(\Phi', \phi, \Phi); \; 1) \; \Phi'= (\beta, \phi); \; 
 \delta \Phi'= (\delta \beta, \phi) + (\beta, \delta \phi); \; 
 (II): \; 0= (\delta \beta, \delta \phi); 
\\
&& 2a) \; \phi =(\Phi', \beta), \; 
 \delta \phi =(\delta \Phi', \beta) + (\Phi', \delta \beta); \; 
 (II): \;  0= (\delta \Phi', \delta \beta). 
\\
&& 2b) \; \phi= (\beta, \Phi); \; 
 \delta \phi= (\delta \beta, \Phi) + (\beta, \delta \Phi); \;   
(II): \; 0= (\delta \beta, \delta \Phi).  
\\
&& 3)\;  \phi=( \Phi', \psi, \Phi); \; 
 \delta \phi=( \delta \Phi', \psi, \Phi) + ( \Phi', \delta \psi, \Phi)
 + ( \Phi', \psi, \delta \Phi); \; (III). 
\end{eqnarray*}
The case $(I00)$ is equivalent to $(I0)$. 
The case $(00I)$ is equivalent to $(0I)$.
Consider the case $(I0I)$. The differential part: 
$(I0I): \; 0=(\delta \Phi', \phi, \delta \Phi); \; (III)$. 
Non-differential part: 
\begin{eqnarray*}
&& (I0I): \; 0=(\delta \Phi', \phi, \delta \Phi); \; 
  1) \; \delta \Phi' =( \beta', \phi); \;(II). 
\\
&& 2a) \; \phi=(\delta \Phi', \beta); \;  
  \delta \phi=(\delta \Phi', \delta \beta); \; 0; \;  
  2b)\; \phi=(\beta', \delta \Phi); \; 
 \delta \phi=(\delta \beta', \delta \Phi); \; 0; 
\\
&& 
3) \; \phi=(\delta \Phi', \gamma, \delta \Phi), \; 
 \delta \phi=(\delta \Phi', \delta \gamma, \delta \Phi); \; 0;  
 \; 4)\;  \delta \Phi=( \phi, \beta); \; (II).  
\end{eqnarray*} 
The case $(II0)$. The differential part: 
$(II0): \; 0=(\delta \Phi', \delta \phi, \Phi); \; (III)$. 
Non-differential part: 
\begin{eqnarray*}
&& (II0): \; 0=(\delta \Phi', \delta \phi, \Phi);  \; 
1) \; \delta \Phi'=(\beta', \delta \phi); \; (II). 
\\
&& 2a) \;  \delta \phi =(\delta \Phi',  \beta); \; (II); \; 
 2b) \; \delta \phi=(\beta',  \Phi); \;  (II); 
\\
&&
 2c) \; \delta \phi =(\delta \Phi', \gamma, \Phi); \; 0; \;  
  3) \;  \Phi=( \delta \phi, \beta);  \;
 \delta \Phi=( \delta \phi, \delta \beta); \; 0. 
\end{eqnarray*}
The case $(0II)$: differential part: 
$(0II): \; 0=( \Phi', \delta \phi, \delta \Phi); \; (III)$.  
Non-differential part: 
\begin{eqnarray*}
&& (0II): \; 0=( \Phi', \delta \phi, \delta \Phi); \; 
 1) \; \Phi'=( \beta',  \delta \phi); \; 
 \delta \Phi'=( \delta \beta',  \delta \phi); \; 0; 
\\
&& 2a) \;  \delta \phi =( \Phi', \beta); \; (II); \; 
 2b) \; \delta \phi =( \beta',  \delta \Phi); \; (II); \; 
2c) \; \delta \phi=( \Phi', \gamma, \delta \Phi); 
\\
&&  0= (\delta \Phi', \gamma, \delta \Phi) + ( \Phi', \delta \gamma, \delta \Phi); \; 
 0=(III); \;
3) \; \delta \Phi =( \delta \phi, \beta); \; (II). 
\end{eqnarray*}
In the case of the sequence ending by $(II)$,
we arrive at an infinite series of differential and orthogonality conditions. 
Namely, with $k \in \Z$ (since a $(II)$
sequence can be seen infinite in both directions),  
we have the sequence 
\begin{eqnarray}
\label{lobuda}
0=(\alpha_k, \beta_k), \quad 
 \delta \alpha_k=( \alpha_{k+1}, \beta_k), \quad    
 \delta \beta_k=(\alpha_k, \beta_{k+1}). 
\end{eqnarray}
Recall, that for each differential relation in this example, 
there exists the coherence relation for corresponding $\mathcal C$-indices. 
Here $\alpha_{k+1}$ and $\beta_{k+1}$ belong to the   
spaces of $\mathcal C$ 
corresponding to the relations \eqref{lobuda}  
via the coherence relations. 
The sequences of differential and orthogonality conditions
derived above 
together with the corresponding coherences relations for indices 
taken for all  
choices of $J$ provide the full set relations defining 
a multiply graded differential algebra.  
\end{example}
The structure of relations for a resulting differential algebra 
as well as the structure of corresponding invariants 
given by Theorem \ref{marcela} 
 depends which position on the tree of differential and 
orthogonality conditions we start with. 
Then the resulting multiple graded differential algebra 
is a reduction of the full graded differential algebra.  
In addition to that, one can also restrict 
the values of indices of gradings for 
vertical and horizontal parts of a complex.  
E.g., one can set $n\ge 0$ instead of $\Z$. 
This would affect coherence conditions for corresponding identities. 
\section{The invariants}
\label{invariants}
In this Section we give a proof of Theorem \ref{marcela}. 
Recall the notations given in Section \ref{setup}. 
The relations \eqref{tuba}  play the role of orthogonality conditions 
 on $\Phi_j$, $1 \le j \le 4$ 
for invariants \eqref{kordo}. 
\begin{proof}
Let $\phi \in \mathcal C$.  
Now let us show that for  
 $\cdot_l$-product completions $\Phi_j$, $1 \le j \le 4$, 
satisfying quite general restriction \eqref{tuba} for the completion,   
we get invariant cohomology classes \eqref{kordo}
 not depending on the choice of  
$\phi$ such that $\delta \phi \in \mathcal I_{nl}(2)$.   
For some fixed choice of $J$, 
 consider the product 
$\left(\epsilon\left(\Phi_1, \overline{d} d \phi \right)  
  \Phi_2, d \phi, \Phi_3, \phi, \Phi_4\right)$. 
We have 
\begin{eqnarray}
\label{pala}
 &&d \left( \epsilon\left(\Phi_1, \overline{d} d \phi \right), 
\Phi_2, d \phi, \Phi_3, \phi, \Phi_4\right) 
 =
\left(\epsilon \left(d \Phi_1, \overline{d} d\phi \right), 
 \Phi_2, d \phi, \Phi_3, \phi, \Phi_4\right) 
\\
&&\quad +
\left(\epsilon \left(\Phi_1, \overline{d} d \phi \right), 
d\Phi_2, d \phi, \Phi_3,  \phi, \Phi_4\right) 
+\left( \epsilon \left(\Phi_1, \overline{d} d \phi \right), 
 \Phi_2, d \phi, d\Phi_3, \phi, \Phi_4\right)  
\nn
\nonumber
&&\quad + 
\left(\epsilon\left(\Phi_1, \overline{d} d \phi \right),
 \Phi_2,d \phi, \Phi_3, d \phi, \Phi_4\right)
+ 
(\epsilon \left(\Phi_1, \overline{d} d \phi \right),
 \Phi_2, d \phi, \Phi_3,  \phi, d \Phi_4)=0,  
\end{eqnarray}
when $d \phi \in \mathcal I_{nl}(2)$, and \eqref{tuba} is satisfied, i.e., 
\begin{equation}
\label{porolon}
\sum_{j=1}^4  d_j.
\left( \epsilon \left( \Phi_1, \overline{d} d \phi \right),  
\Phi_2, d \phi, \Phi_3, \phi, \Phi_4 \right)=0,  
\end{equation}
so that the differential of the total completion is zero. 
Then \eqref{pala} defines the cohomology class  
$\left[\epsilon\left(\Phi_1, \overline{d} d \phi \right),
 \Phi_2, d \phi, \Phi_3, \phi, \Phi_4)\right]$. 

Now, let us check that for completions $\Phi_j$, $1 \le j \le 4$, 
satisfying \eqref{tuba},  
the cohomology classes \eqref{kordo} do not depend
 on the choice of $\phi$.   
 Indeed, let us substitute 
$\phi$ by $\phi + \eta$, with $\delta \eta \in \mathcal I_{nl}(2)$,  
and satisfying conditions \eqref{tuba}. 
Then,   
\begin{eqnarray*}
&& \left(\epsilon\left(\Phi_1, \overline{d} d (\phi+\eta)\right), \Phi_2,  
d (\phi+ \eta), \Phi_3, \phi+\eta, \Phi_4\right) 
\nn
&&= 
\left(\epsilon\left(\Phi_1, \overline{d} d \phi\right), \Phi_2, 
d \phi, \Phi_3, \phi, \Phi_4\right) 
 +\left(\epsilon\left(\Phi_1, \overline{d} d \phi\right), \Phi_2, 
d \phi, \Phi_3, \eta, \Phi_4\right) 
\nn
&&+\left(\epsilon\left(\Phi_1, \overline{d} d \phi\right), \Phi_2, 
d  \eta, \Phi_3, \phi, \Phi_4\right)  
+
\left(\epsilon\left(\Phi_1, \overline{d} d \phi\right), \Phi_2, 
d  \eta, \Phi_3, \eta, \Phi_4\right) 
\nn
&&+
\left(\epsilon\left(\Phi_1, \overline{d} d\eta\right), \Phi_2, 
d \phi, \Phi_3, \phi, \Phi_4\right)
+
\left(\epsilon\left(\Phi_1, \overline{d} d\eta\right), \Phi_2, 
d \phi, \Phi_3, \eta, \Phi_4\right)
\nn
&&+
\left(\epsilon\left(\Phi_1, \overline{d} d \eta\right), \Phi_2, 
d \eta, \Phi_3, \phi, \Phi_4\right) 
+
 \left(\epsilon\left(\Phi_1, \overline{d} d \eta\right), \Phi_2, 
d \eta, \Phi_3, \eta, \Phi_4\right). 
\end{eqnarray*}
Note that 
\begin{eqnarray*}  
 && \left(\epsilon\left(d \Phi_1, \overline{d}  \phi\right), 
\Phi_2,  d \phi, \Phi_3, \eta, \Phi_4\right)
 + \left(\epsilon\left(\Phi_1, \overline{d} d \phi\right),
  \Phi_2,  d \phi, \Phi_3, \eta,   \Phi_4\right)
\nn
&&
 + \left(\epsilon \left(\Phi_1, \overline{d}  \phi\right),
  d \Phi_2,  d \phi,  \Phi_3, \eta,  \Phi_4\right) 
  + \left(\epsilon \left(\Phi_1, \overline{d} \phi\right),
 \Phi_2,  d \phi,  d \Phi_3,  \eta,  \Phi_4\right)    
\nn
&& 
+ \left(\epsilon \left(\Phi_1, \overline{d} \phi\right), 
\Phi_2,  d \phi, \Phi_3,  d \eta,  \Phi_4\right)  
+ \left(\epsilon\left(\Phi_1, \overline{d}  \phi\right),
 \Phi_2,  d \phi, \Phi_3,  \eta, d \Phi_4\right) 
\nn
&& =d \left(\epsilon\left(\Phi_1, \overline{d}  \phi\right), 
 \Phi_2,  d \phi, \Phi_3, \eta, \Phi_4\right) 
=\left(\epsilon\left(\Phi_1, \overline{d} d \phi\right), 
  \Phi_2,  d \phi, \Phi_3, \eta,   \Phi_4\right), 
\end{eqnarray*}
and similarly  
\begin{eqnarray*}
&& d \left(\epsilon \left( \Phi_1, \overline{d}  \phi \right), \Phi_2, 
d  \eta, \Phi_3, \phi, \Phi_4 \right) 
=
\left(\epsilon\left(\Phi_1, \overline{d} d \phi\right), \Phi_2, 
d  \eta, \Phi_3, \phi, \Phi_4 \right),   
\nn
&& d \left(\epsilon \left(\Phi_1, \overline{d}  \phi \right), \Phi_2,  
d  \eta, \Phi_3, \eta, \Phi_4 \right) 
=\left(\epsilon \left( \Phi_1, \overline{d} d \phi \right), \Phi_2, 
d  \eta, \Phi_3, \eta, \Phi_4\right),  
\nn
&& d \left( \epsilon \left( \Phi_1, \overline{d} \eta\right), \Phi_2, 
d \phi, \Phi_3, \phi, \Phi_4\right)
=\left(\epsilon\left(\Phi_1, \overline{d} d\eta\right), \Phi_2, 
d \phi, \Phi_3, \phi, \Phi_4\right), 
\nn 
&& d \left( \epsilon \left( \Phi_1, \overline{d} \eta \right), \Phi_2, 
d \phi, \Phi_3, \eta, \Phi_4 \right)
=
\left( \epsilon \left( \Phi_1, \overline{d} d\eta \right), \Phi_2,  
d \phi, \Phi_3, \eta, \Phi_4 \right), 
\nn
&& d \left(\epsilon \left(\Phi_1, \overline{d}  \eta\right), \Phi_2, 
d \eta, \Phi_3, \phi, \Phi_4\right) 
=\left(\epsilon \left(\Phi_1, \overline{d} d \eta\right), \Phi_2, 
d \eta, \Phi_3, \phi, \Phi_4\right), 
\nn
 &&
 d \left( \epsilon \left(\Phi_1, \overline{d} \eta \right), \Phi_2,   
 d \eta, \Phi_3, \eta, \Phi_4 \right)  
= \left( \epsilon \left( \Phi_1, \overline{d} d \eta \right),
 \Phi_2, 
d \eta, \Phi_3, \eta, \Phi_4  \right), 
\end{eqnarray*}
due to \eqref{porolon} and $d\phi$, $d\eta \in \mathcal I_{nl}(2)$. 
Therefore, 
\begin{eqnarray*}
&&\left(\epsilon\left(\Phi_1, \overline{d} d (\phi+\eta)\right), \Phi_2, 
d (\phi+ \eta), \Phi_3, \phi+\eta, \Phi_4\right)  
\nn
&=& \left(\epsilon \left(\Phi_1, \overline{d} d \phi\right), \Phi_2, 
d \phi, \Phi_3, \phi, \Phi_4\right) 
+ d \left( \left(
\epsilon \left(\Phi_1, \overline{d}  \phi\right), \Phi_2, 
d \phi, \Phi_3, \eta, \Phi_4 \right)  \right. 
\nn
&+& \left(\epsilon \left(\Phi_1, \overline{d} \phi \right), \Phi_2, 
d  \eta, \Phi_3, \phi, \Phi_4\right)  
+ \left(\epsilon \left(\Phi_1, \overline{d}  \phi\right), \Phi_2, 
d  \eta, \Phi_3, \eta, \Phi_4\right)  
\nn
&+&
\left(\epsilon \left(\Phi_1, \overline{d} \eta\right), \Phi_2, 
d \phi, \Phi_3, \phi, \Phi_4\right)
+
\left(\epsilon \left(\Phi_1, \overline{d} \eta\right), \Phi_2, 
d \phi, \Phi_3, \eta, \Phi_4\right)
\nn
&+& \left.
\left(\epsilon \left(\Phi_1, \overline{d} \eta\right), \Phi_2, 
d \eta, \Phi_3, \phi, \Phi_4\right) 
+
 \left(\epsilon \left(\Phi_1, \overline{d}  \eta\right), \Phi_2, 
d \eta, \Phi_3, \eta, \Phi_4 \right) \right).  
\end{eqnarray*} 
Thus, the cohomology class is preserved under the transformations  
$\phi\to \phi+ \eta$. 
In the opposite choice of $d$, we use the symmetry of $\delta$ and $\Delta$
for vertical chain-cochain complexes with 
 differentials $\Delta$ commuting with differentials $\delta$.  
The considerations above extend to the symmetrized expression \eqref{kordo}. 
\end{proof}
\section*{Acknowledgment}
The second author is supported by 
the Institute of Mathematics, Academy of Sciences of the Czech Republic (RVO 67985840).  


\begin{thebibliography}{99}
\bibitem
{BG}  Bazaikin, Ya. V.,  Galaev, A. S.  Losik classes for codimension one foliations, 
 J. Inst. Math. Jussieu 21 (2022), no. 4, 1391--1419. 
\bibitem
{BGG} Ya. V. Bazaikin, A. S. Galaev, and P. Gumenyuk. 
Non-diffeomorphic Reeb foliations and modified Godbillon-Vey class, 
Math. Z. (2022), no. 2, 1335--1349.


\bibitem
{BGG} Ya. V. Bazaikin, A. S. Galaev, N. I. Zukova. 
 Chaos in Cartan foliations. Chaos 30 (2020), no. 10, 103116. 

\bibitem{CM} M. Crainic and I. Moerdijk. 
 ${\rm \check C}$ech-De Rham theory for leaf spaces of foliations. Math.
Ann. 328 (2004), no. 1--2, 59--85.


\bibitem 
{BZF} 
E. Frenkel, D. Ben-Zvi. Vertex algebras and algebraic curves. Mathematical Surveys and Monographs,
 88. American Mathematical Society, Providence, RI, 2001. xii+348 pp.
 

\bibitem
{FMS} Ph. Francesco, P. Mathieu, and D. Senechal. Conformal Field Theory. 
 Graduate Texts in Contemporary Physics. 1997. 

\bibitem{Galaev} 
A. S. Galaev.  
Comparison of approaches to characteristic classes of foliations. arXiv:1709.05888. 


\bibitem{Ghys} E. Ghys.     
L'invariant de Godbillon-Vey. Seminaire Bourbaki, 
41--eme annee, n 706, S. M. F.
Asterisque 177--178 (1989). 


\bibitem{G}  
B. Gui. Convergence of sewing conformal blocks. Comm. Contemp. Math. 
Vol. 26, No. 03, 2350007 (2024). 

\bibitem{Ko} 
D. Kotschick. Godbillon-Vey invariants for families of foliations. Symplectic and contact 
topology: interactions and perspectives (Toronto, ON/Montreal, QC, 2001), 
131--144, Fields Inst. Commun., 35, Amer. Math. Soc., Providence, RI, 2003. 


\bibitem
{Huang} Y.-Zh. Huang. 
A cohomology theory of grading-restricted vertex algebras. 
Comm. Math. Phys. 327 (2014), no. 1, 279--307. 

\bibitem{Losik} M. V. Losik.   
On some generalization of a manifold and its characteristic classes (Russian), Funcional. Anal. i 
Prilozhen.  
\textbf{24}(1990), no 1, 29-37 ;  English translation in Functional Anal. Appl. \textbf {24} (1990), 26--32. 


\bibitem{sv2} M. V. Saveliev,  A. M. Vershik. 
Continuum analogues of contragredient Lie algebras.
Commun. Math. Phys. 126, 367, 1989. 


\bibitem
{TUY} A. Tsuchiya, K. Ueno, and Y. Yamada, Y. Conformal field
theory on universal family of stable curves with gauge symmetries,
Adv. Stud. Pure. Math. \textbf{19} (1989), 459--566.

\bibitem
{Yamada} A. Yamada. Precise variational formulas for abelian
differentials. Kodai Math.J. \textbf{3} (1980), 114--143.


\end{thebibliography}
\end{document}